\documentclass[a4paper,12pt]{article}
\usepackage[T2A]{fontenc}
\usepackage[cp1251]{inputenc}
\usepackage[dvips]{graphicx}
\usepackage{amssymb,amsmath}
\begin{document}

\title{INTEGRABLE ISOTROPIC GEOMETRICAL FLOWS AND HEISENBERG FERROMAGNETS}
\author{N.S.Serikbaev, Zh.M.Bitibaeva, K.K.Yerzhanov, R.Myrzakulov\footnote{cnlpmyra1954@yahoo.com,
cnlpmyra@mail.ru}}

\maketitle

{\it Department of   General and Theoretical  Physics,
 Eurasian National University, Astana, 010008, Kazakhstan}
\begin{abstract}
Geometrical flows (GF) play an important role in modern mathematics and
physics. In this letter we have considered some integrable
isotropic GF -- Ricci flows (RF) and mean curvature flows (MCF) -- which are related with  integrable
Heisenberg ferromagnets. In 2+1 dimensions, these GF have
a singularity at $t=t_{0}$.
\end{abstract}

\tableofcontents


\section{Introduction}

Geometric flow (GF) is the gradient flow associated with a
functional on a manifold which has a geometric interpretation,
usually associated with some extrinsic or intrinsic curvature. A
GF is also called a geometric evolution equation. They can be
interpreted as flows on a moduli space (for intrinsic flows) or a
parameter space (for extrinsic flows). These are of fundamental
interest in the calculus of variations, and include several famous
problems and theories. Particularly interesting are their critical
points. GF play an important role in mathematics and physics.

In this note we will consider  two examples of GF namely the Ricci
flow (RF) and the mean curvature flow (MCF) related with the
integrable Heisenberg ferromagnets (HF) in 1+1  and 2+1
dimensions. In particular we explore integrable reductions of the
following (2+1)-dimensional GF:
$$
{\bf r}_t=M{\bf n}-{\bf \xi} +u{\bf r}_x,\eqno(1.1a)
$$
$$
u_x=-{\bf r}_x\cdot({\bf r}_{xx}\wedge{\bf r}_{xy})\eqno(1.1b)
$$
and
$$
{\bf r}_{xt}=M{\bf n}-{\bf \eta} +u_{x}{\bf r}_{xy}+u_{y}{\bf
r}_{xx},\eqno(1.2a)
$$
$$
u_{xx}-\alpha^2u_{yy}=-2\alpha^2{\bf r}_x\cdot({\bf
r}_{xx}\wedge{\bf r}_{xy}).\eqno(1.2b)
$$
Here $M(x,y,t)$ and $u(x,y,t)$ are scalar (real) functions. If
$M=H$, then these GF can be considered as the examples of the
(2+1)-dimensional MCF. The particular cases of these GF are the
M-I flow and the Ishimori flow which as MCF are integrable.

\section{Mean curvature flows}
MCF is an example of a GF of hypersurfaces in a Riemannian
manifold (for example, smooth surfaces in 3-dimensional Euclidean
space).
\subsection{Isotropic MCF}
The isotropic MCF reads as (see, e.g. [1])
$$
{\bf r}_{t}=H{\bf n}-{\bf \xi}, \eqno(2.1)
$$
where ${\bf r}=(r_1, r_2, r_3), \quad {\bf n}=(n_1, n_2, n_3)
\quad {\bf \xi}=(\xi_{1}, \xi_{2}, \xi_{3})$ and  $H$ is a mean
curvature. Additionally in this note we assume that
$$
{\bf r}_{x}^{2}=1. \eqno(2.2)
$$
In this note we use also the following form of MCF
$$
{\bf r}_{xt}=H^{\prime}{\bf n}+ \eta, \eqno(2.3)
$$
where $ \eta=(\eta_{1}, \eta_{2}, \eta_{3})$.
\subsection{Anisotropic MCF}
There exist also the
anisotropic MCF which can be written as
$$
{\bf r}_{t}=H{\bf n}-{\bf \xi} +{\bf V},\eqno(2.4a)
$$
$$
{\bf V}_{x}={\bf r}_{x}\wedge J{\bf r}_{x}. \eqno(2.4b)
$$
This equation can be rewritten in the following equivalent form
$${\bf r}_{xt}=H^{\prime}{\bf n}+{\bf \eta}+{\bf r}_{x}\wedge J{\bf r}_{x}, \eqno(2.5)
$$
 where $J=diag(J_{1}, J_{2}, J_{3})$, ${\bf V}=(V_1, V_2, V_3)$.
 Reference [4] presented examples of anisotropic MCF and RF.
\section{Ricci flow}
 \subsection{Isotropic RF}
 The RF is an intrinsic geometric flow—a process which deforms the
metric of a Riemannian manifold—in this case in a manner formally
analogous to the diffusion of heat, thereby smoothing out
irregularities in the metric. It plays an important role in the
proof of the Poincare conjecture.

 Given a Riemannian manifold with
metric tensor $g_{ij}$, we can compute the Ricci tensor $R_{ij}$,
which collects averages of sectional curvatures into a kind of
"trace" of the Riemann curvature tensor. If we consider the metric
tensor (and the associated Ricci tensor) to be functions of a
variable which is usually called "time" (but which may have
nothing to do with any physical time), then the RF may be defined
by the geometric evolution equation
$$
g_{ijt}=-2R_{ij}.\eqno(3.1a)
$$
At the same time, the equation of the normalized RF reads
$$
g_{ijt}=-2R_{ij}+\frac{2}{n}<R>g_{ij},\eqno(3.1b)
$$
where $<R>$ is the average (mean) of the scalar curvature (which
is obtained from the Ricci tensor by taking the trace) and $n$ is
the dimension of the manifold. Equation (3.2) preserves the
volume of the metric.
\subsection{Anisotropic RF}
In the anisotropic case, Eqs. (3.1) take the form
$$
g_{ijt}=-2R_{ij}+A_{ij}\eqno(3.2a)
$$
and
$$
g_{ijt}=-2R_{ij}+\frac{2}{n}<R>g_{ij}+A_{ij},\eqno(3.2b)
$$
respectively.
\section{HF flow} Consider the following HF
$$
{\bf S}_{y}={\bf S}\wedge{\bf S}_{xx}, \eqno(4.1)
$$
where ${\bf S}=(S_{1},S_{2},S_{3}), \quad {\bf S}^{2}=1$. If we
assume
$$
{\bf S}={\bf r}_{x},  \eqno(4.2)
$$
then the equation (4.1) takes the form [5]
$$
{\bf r}_{yx}=({\bf r}_{x}\wedge{\bf r}_{xx})_x \eqno(4.3a)
$$
or
$$
{\bf r}_{y}={\bf r}_{x}\wedge{\bf r}_{xx}. \eqno(4.3b)
$$
The corresponding Lax representation is given by
$$
\Phi_x=U\Phi,\quad \Phi_y=V\Phi, \eqno(4.4)
$$
where
$$
U=\frac{\lambda}{2i}r_{x},\quad
V=\frac{i\lambda^2}{2}r_x+\frac{\lambda}{2}r_{xx}r_{x},\quad
r={\bf r}\cdot \mathbf{\sigma},\quad {\bf
\sigma}=(\sigma_1,\sigma_2,\sigma_3).
$$
For the equation (4.3) we get
$$
E=1, \quad F=0, \quad G={\bf r}^{2}_y.\eqno(4.5)
$$
Here $E,F,G$ and $L,M,N$ are the coefficients of the fundamental
forms of the surface
$$
I=d{\bf r}^2=g_{ij}dx^i dx^j=Edx^2+2Fdxdy+Gdy^2,\eqno(4.6a)
$$
$$
 II=d{\bf r}\cdot{\bf n}=b_{ij}dx^i dx^j=Ldx^2+2Mdxdy+Ndy^2.\eqno(4.6b)
$$
In our case we have
$$
R_{ij}= \frac{1}{2}Rg_{ij},\eqno(4.7a)
$$
$$
R= \frac{G^2_x-2GG_{xx}}{2G^2}, \eqno(4.7b)
$$
$$
K= \frac{R}{2}.\eqno(4.7c)
$$
The modified  RF related with the HF (4.1) can be written in the
following form
$$ g_{ijt}=-2R_{ij}+F_{ij}.\eqno(4.8)
$$
Now we present the MCF related with the HF equation (4.1). To do
it, we consider the surfaces in $R^3$ associated to two parameters
$x$ and $y$ and the renormalization group time $t$. It is
convenient, where appropriate, to think of the surface as a graph of
a function $r_3=\varphi(r_1,(x,y;t),r_2(x,y;t);t)$ that evolves in
time. In our case
$$
r_{1y}=r_{2x}r_{3xx}-r_{xx}r_{3x},\eqno(4.9a)
$$ $$
r_{2y}=r_{3x}r_{1xx}-r_{3xx}r_{1x},\eqno(4.9b)
$$
$$
r_{3y}=r_{1x}r_{2xx}-r_{1xx}r_{2x}.\eqno(4.9c)
$$ In this system we must use the following expressions
$$
r_{3y}=\varphi_{r_1}r_{1y}+\varphi_{r_2}r_{2y},\eqno(4.10a)
$$
$$
r_{3x}=\varphi_{r_1}r_{1x}+\varphi_{r_2}r_{2x},  \eqno(4.10b)
$$
$$
r_{3xx}=\varphi_{r_1r_1}r^{2}_{1x}+2\varphi_{r_1}\varphi_{r_2}r_{1x}r_{2x}+\varphi_{r_2r_2}r_{2x}^2,\eqno(4.10c)
$$
From these and (2.2) we can define two functions $r_{1x},
r_{1y}$ as
$$
r_{1x}=\frac{-\varphi_{r_1}\varphi_{r_2}r_{2x}\pm\sqrt{1+\varphi_{r_1}^{2}-(1+\varphi_{r_1}^{2}+
\varphi_{r_2}^2)r^{2}_{2x}}}{1+\varphi_{r_1}^{2}},\eqno(4.11a)
$$
$$
r_{1y}=\frac{r_{1x}r_{2xx}-r_{1xx}r_{2x}-\varphi_{r_2}r_{2y}}{\varphi_{r_1}}.\eqno(4.11b)
$$
In this notation the mean curvature of the surface can be
written as [1]
$$
H=\frac{(1+(\varphi_{r_2})^2)\varphi_{r_1r_1}+(1+(\varphi_{r_1})^2)\varphi{r_2r_2}-2\varphi_{r_1}\varphi_{r_2}
\varphi_{r_2r_2}}{(\sqrt{1+(\varphi{r_1})^2+(\varphi_{r_2})^2})^3}.\eqno(4.12)
$$
The inward unit normal vector is given by
$$
{\bf n}=\frac{1}{\sqrt{1+(\varphi_{r_1})^2+(\varphi_{r_2})^2}}
(-\varphi_{r_1},-\varphi_{r_2},1).\eqno(4.13)
$$
So for the MCF we have the following equation [1]
$$
\varphi_t=\frac{(1+(\varphi_{r_2})^2)\varphi_{r_1r_1}+(1+\varphi_{r_1}^2)\varphi_{r_2r_2}^2-2\varphi_{r_1}\varphi_{r_2}
\varphi_{r_2r_2}}{1+(\varphi_{r_1})^2+(\varphi_{r_2})^2}+
$$
$$
\xi_{1}\varphi_{r_1}+\xi_{2}\varphi_{r_2}-\xi_{3}.\eqno(4.14)
$$
For the first and second fundamental forms of the
two-dimensional surface and for the mean curvature we have the
following system of equations [1]
$$
g_{ijt}=-2HK_{ij},\eqno(4.15a)
$$
$$
g^{ij}_t=2HK^{ij}.\eqno{4.15b}
$$
Hence we get
$$
(\ln\sqrt{g})_{t}=-H^2,\eqno(4.16a)
$$
$$
K_{ijt}=g^{lm}\nabla_l\nabla_mK_{ij}-2H(K^2)_{ij}+(TrK^2)K_{ij},\eqno(4.16b)
$$
$$
\frac{\partial H}{\partial
t}=g^{ij}\nabla_i\nabla_jH+(TrK^2)H.\eqno(4.16c)
$$

\section{M-I flow as MCF}
Now we consider the following Myrzakulov I equation
(abbreviated as the M-I equation)
$$
{\bf S}_t=({\bf S}\wedge{\bf S}_y+u{\bf S})_x,\eqno(5.1a)
$$
$$
u_x=-{\bf S}\cdot({\bf S}_x\wedge{\bf S}_y).\eqno(5.1b)
$$
After the identification (4.2) this system takes the form
$$
{\bf r}_t={\bf r}_x\wedge{\bf r}_{xy}+u{\bf r}_x,\eqno(5.2a)
$$
$$
u_x=-{\bf r}_x\cdot({\bf r}_{xx}\wedge{\bf r}_{xy}).\eqno(5.2b)
$$
This is M-I flow and is known to be integrable. The corresponding Lax
representation is given by
$$
\Phi_{x}=U\Phi,\eqno(5.3a)
$$
$$
\Phi_{t}=\lambda\Phi_{y}+V\Phi,\eqno(5.3b)
$$
where
$$
U=\frac{i\lambda}{2}r_{x},\eqno(5.4a)
$$
$$
V=\frac{\lambda}{4}([r_{x},r_{xy}]+2iur_{x}).\eqno(5.4b)
$$
Note that the M-I flow (5.2) can be written in the following form
$$
{\bf r}_t=\left(\frac{MF}{\sqrt{g}}+u\right){\bf
r}_x-\frac{M}{\sqrt{g}}{\bf r}_y+\frac{G_x}{2\sqrt{g}}{\bf n}+{\bf
c},\eqno(5.5a)
$$
$$
u_x=\frac{LG_x-2MF_x}{2\sqrt{g}}.\eqno(5.5b)
$$
Here ${\bf c}={\bf c}(y,t)$ which we set to ${\bf c}=0$. For the
M-I flow we have
$$
g_{11t}=E_{t}=0,\eqno(5.6a)
$$
$$
g_{12t}=F_{t}=\frac{(FM-N)F_x}{\sqrt{g}}-L_y\sqrt{g}+uF_x+u_y,\eqno(5.6b)
$$
$$
g_{22t}=G_{t}=\frac{(FM-N)G_x}{\sqrt{g}}-2M_y\sqrt{g}+uG_x+2Fu_y.\eqno(5.6c)
$$
So that for  $g=\det(g_{ij})$ we get
$$
g_t=G_t-2FF_t=\frac{g_x}{\sqrt{g}}(FM-N)-2\sqrt{g}(M_y-FL_y)+ug_x. \eqno(5.7)
$$
Equations (5.6)-(5.7) are integrable.  Note that in our case
we have the following formula
$$
R_{ij}=\frac{1}{2}Rg_{ij}\eqno(5.8)
$$
and
$$
R=\frac{1}{2g^2}\left(4FF_xF_y-2F_xG_y-2FF_xG_x+G^2_x-4F^2F_{xy}+4GF_{xy}+2F^2G_{xx}-2GG_{xx}\right).\eqno(5.9)
$$
As in the previous case it is convenient, where appropriate, to
think of the surface as the graph of a function
$r_3=\varphi(r_1(x,y;t),r_2(x,y;t);t)$ that evolves in time. First
let us rewrite the system (5.2) in component form. We have
$$
r_{1t}=r_{2x}r_{3xy}-r_{2xy}r_{3x}+ur_{1x},\eqno(5.10a)
$$
$$
r_{2t}=r_{3x}r_{1xy}-r_{3xy}r_{3x}+ur_{2x},\eqno(5.10b)
$$
$$
r_{3t}=r_{1x}r_{2xy}-r_{1xy}r_{2x}+ur_{3x},\eqno(5.10c)
$$
$$
u_{x}=-r_{1xx}(r_{2xx}r_{3xy}-r_{2xy}r_{3xx})-r_{2xx}(r_{3xx}r_{1xy}-$$
$$r_{3xy}r_{1xx})-r_{3xx}(r_{1xx}r_{2xy}-r_{1xy}r_{2xx}).\eqno(5.10d)
$$
As
$$ r_{3y}=\varphi_{r_1}r_{2y}+\varphi_{r_2}r_{2y},\eqno(5.11a)
$$
$$
r_{3x}=\varphi_{r_1}r_{1x}+\varphi_{r_2}r_{2x},  \eqno(5.11b)
$$
$$
r_{3xx}=\varphi_{r_1r_1}r^{2}_{1x}+2\varphi_{r_1}\varphi_{r_2}r_{1x}r_{2x}+\varphi_{r_2r_2}r_{2x}^{2},\eqno(5.11c)
$$
for the  functions $r_{1x}, r_{1y}$ we get
$$
r_{1x}=\frac{-\varphi_{r_1}\varphi_{r_2}r_{2x}\pm\sqrt{1+\varphi_{r_1}^{2}-(1+\varphi_{r_1}^{2}+\varphi_{r_2}^{2})r^{2}_{2x}}}{1+\varphi_{r_1}^{2}},\eqno(5.12a)
$$
$$
r_{1y}=\frac{r_{1x}r_{2xx}-r_{1xx}r_{2x}-\varphi_{r_2}r_{2y}}{\varphi_{r_1}}.\eqno(5.12b)
$$
So for the mean curvature of the surface we have the following
formula
$$
H=\frac{(1+(\varphi_{r_2})^2)\varphi_{r_1}^2+(1+\varphi^2_{r_1})\varphi_{r_2}^2-2\varphi_{r_1}\varphi{r_2}
\varphi_{r_1r_2}}{(\sqrt{1+\varphi^2_{r_1}+\varphi^2_{r_2}})^3}.\eqno(5.13)
$$
On the other hand the inward unit normal vector is defined by the
formula [1]
$$
{\bf n}=\frac{1}{\sqrt{1+\varphi^2_{r_1}+\varphi^2_{r_2}}}(-\varphi_{r_1},-\varphi_{r_2},1).\eqno(5.14)
$$
Now we are ready to find the MCF. As result we obtain
$$
\varphi_t=\frac{(1+(\varphi_{r_2})^2)\varphi_{r_1r_1}+(1+(\varphi_{r_1})^2)\varphi_{r_2r_2}^2-2\varphi_{r_1}\varphi_{r_2}
\varphi_{r_2r_2}}{1+(\varphi_{r_1})^2+(\varphi_{r_2})^2}+
$$
$$
\xi_{1}\varphi_{r_1}+\xi_{2}\varphi_{r_2}-\xi_{3}.\eqno(5.15)
$$
Finally we present the following two equations [1]
$$
(\ln\sqrt{g})_{tt}=-2Hg^{ij}\nabla_i\nabla_jH-2(TrK^2)H^2\eqno(5.16)
$$
and
$$
(\sqrt{g})_{tt}=[-2Hg^{ij}\nabla_i\nabla_jH-2(TrK^2)H^2+H^4]\sqrt{g}.\eqno(5.17)
$$
\section{M-I flow as RF}

Consider a Riemannian manifold $M$ of dimension $3$ with a local
coordinate system $X^\mu$ and metric $G_{\mu\nu}(X)$ so that its
first fundamental form is
$$
ds^2_{\it M}=G_{\mu\nu}(X)dX^\mu dX^\nu \eqno (6.1a)
$$
or
$$
ds^2=G_{11}dx^2+G_{22}dy^2+G_{33}dt^2+2G_{12}dxdy+2G_{13}dxdt+2G_{23}dydt.\eqno(6.1b)
$$
Here $X^1=x$, $X^2=y$, $X^3=t$ and
$$ G_{11}=\mathbf{r}^{2}_{x},\quad
G_{22}=\mathbf{r}^{2}_{y},\quad G_{33}=\mathbf{r}^{2}_{t},
$$
$$ G_{12}=G_{21}=\mathbf{r}_{x}\cdot\mathbf{r}_{y}, \quad
G_{13}=G_{31}=\mathbf{r}_{x}\cdot\mathbf{r}_{t},\quad
 G_{23}=G_{32}=\mathbf{r}_{y}\cdot\mathbf{r}_{t}.\eqno(6.2)
$$
For the M-I flow the metric takes the form
$$ G_{11}=\mathbf{r}^{2}_{x}=1,\quad G_{22}=\mathbf{r}^{2}_{y},\quad
 G_{33}=\mathbf{r}^{2}_{xy}+u,
$$
$$ G_{12}=G_{21}=\mathbf{r}_{x}\cdot\mathbf{r}_{y}, \quad
G_{13}=G_{31}=u,\quad
 G_{23}=G_{32}=u\mathbf{r}_{y}\cdot\mathbf{r}_{y}-\mathbf{r}_{xy}\cdot(\mathbf{r}_{x}\wedge\mathbf{r}_{y}).\eqno(6.3)
$$
In this case $G=\det(G_{ij})$ is given by
$$G=u^2[\mathbf{r}_{xy}^2-\mathbf{r}_{y}^2]+u[\mathbf{r}_{y}^2-(\mathbf{r}_{x}\cdot\mathbf{r}_{y})^2]+[\mathbf{r}_{y}^2-(\mathbf{r}_{x}\cdot\mathbf{r}_{y})^2]
\mathbf{r}_{xy}^2.\eqno(6.4)
$$
For the M-I flow the corresponding
modified RF reads as
$$
G_{ijt}=-2R_{ij}+F_{ij}. \eqno(6.5)
$$
This modified RF is integrable. This follows from the integrability of the original equation (5.2).
\section{Conclusion}
In this letter we have considered some integrable and non
integrable  GF related with some integrable HF in 1+1 and 2+1.

Also we would like to note that the flows considered in this letter
have a singularity with respect to $t$. This  is related with the
following fact: The spectral parameter $\lambda$ in the Lax
representation (5.3) obeys the following nonlinear equation
$$
\lambda_{t}=\lambda\lambda_{y}.\eqno(7.1)
$$
This equation has the following solution
$$
\lambda=\frac{a+y}{t_0-t}.\eqno(7.2)
$$
It has a singularity at $t=t_{0}$. This means that RF and MCF,
which are related with the M-I flow, also have a singularity at
this point. Also we note that references [2]-[3] considered GF
related to the Ishimori equation.  Finally we would like to
present the following systems [3]
$$
g_{ijt}=-2R_{ij}+ug_{ij}, \eqno(7.3a)
$$
$$
u_{tt}=\triangle (u+kR)\eqno(7.3b) $$ and

$$
g_{ijt}=-2R_{ij}+ug_{ij}, \eqno(7.4a)
$$
$$
u_{t}=\triangle (u+kR).\eqno(7.4b) $$ More general forms of these
systems look like [3]
$$
g_{ijt}=-2R_{ij}+(\beta u+\alpha R)g_{ij}, \eqno(7.5a)
$$
$$
u_{tt}=\triangle (u+kR)\eqno(7.5b) $$ and

$$
g_{ijt}=-2R_{ij}+(\beta u+\alpha R)g_{ij}, \eqno(7.6a)
$$
$$
u_{t}=\triangle (u+kR),\eqno(7.6b) $$ respectively.

\end{document}